\documentclass[11pt]{amsart}
\newtheorem*{fit}{The First Incompleteness Theorem}
\newtheorem*{sit}{The Second Incompleteness Theorem}
\newtheorem*{ct}{The Completeness Theorem}
\newtheorem*{theorem}{Theorem}
\theoremstyle{definition}
\newtheorem*{definition}{Definition}
\begin{document}
\noindent
{\small Topology Atlas Invited Contributions \textbf{9} no.~3 (2004) 6 pp.}
\vskip 0.25in
\title[Hilbert's problems and the foundations of mathematics]{Hilbert's 
First and Second Problems and the foundations of mathematics}
\author{Peter Nyikos}
\address{University of South Carolina, Columbia, SC 29208 USA}
\thanks{This article is adapted from a one-hour colloquium lecture given
at the University of Auckland in May, 2000, just three months before the 
100th anniversary of Hilbert's lecture.}
\maketitle

In \oldstylenums{1900}, David Hilbert gave a seminal lecture in which he
spoke about a list of unsolved problems in mathematics that he deemed to
be of outstanding importance.  The first of these was Cantor's continuum
problem, which has to do with infinite numbers with which Cantor
revolutionised set theory.  The smallest infinite number, $\aleph_0$,
`aleph-nought,' gives the number of positive whole numbers.  A set is of
this cardinality if it is possible to list its members in an arrangement
such that each one is encountered after a finite number (however large) of
steps.  Cantor's revolutionary discovery was that the points on a line
cannot be so listed, and so the number of points on a line is a strictly
higher infinite number ($\mathfrak{c}$, `the cardinality of the
continuum') than $\aleph_0$.  Hilbert's First Problem asks whether any
infinite subset of the real line is of one of these two cardinalities.  
The axiom that this is indeed the case is known as the Continuum
Hypothesis (\textsf{CH}).

This problem had unexpected connections with Hilbert's Second Problem (and
even with the Tenth, see the article by \textsc{M. Davis} and the
comments on the book edited by \textsc{F. Browder}). The Second Problem
asked for a proof of the consistency of the foundations of mathematics.  
Some of the flavor of the urgency of that problem is provided by the
following passage from an article by \textsc{S.G. Simpson} in the same 
volume of JSL as the article by \textsc{P. Maddy}:

\begin{quote}
We must remember that in Hilbert's time, all mathematicians were excited
about the foundations of mathematics.  Intense controversy centered around
the problem of the legitimacy of abstract objects. Weierstrass had greatly
clarified the role of the infinite in calculus.  Cantor's set theory
promised to raise mathematics to new heights of generality, clarity and
rigor.  But Frege's attempt to base mathematics on a general theory of
properties led to an embarrassing contradiction. Great mathematicians such
as Kronecker, Poincar\'e, and Brouwer challenged the validity of all
infinitistic reasoning.  Hilbert vowed to defend the Cantorian paradise.  
The fires of controversy were fueled by revolutionary developments in
mathematical physics. There was a stormy climate of debate and criticism.
The contrast with today's climate of intellectual exhaustion and
compartmentalization could not be more striking.

\dots Actually, Hilbert saw the issue as having supramathematical
significance.  Mathematics is not only the most logical and rigorous of
the sciences but also the most spectacular example of the power of
``unaided'' human reason.  If mathematics fails, then so does the human
spirit.  I was deeply moved by the following passage [13, pp.~370--371]:
``The definitive clarification of the nature of the infinite has become
necessary, not merely for the special interests of the individual sciences
but for the honor of human understanding itself.''
\end{quote}

Hilbert was already aware, at the time of his \oldstylenums{1900} lecture,
of some connection between the provability of the consistency of a
mathematical theory and the decidability of statements by the axioms of
the theory. But it was Kurt G\"odel who showed the true nature of this
connection in the process of showing that Hilbert's Second Problem has a
negative solution:

\begin{fit}
Every recursively axiomatizable theory rich enough to include the Peano 
Axioms contains statements whose truth cannot be decided within the 
theory.
In particular, Peano Arithmetic itself can be used to formulate true 
statements about the natural numbers that are not provable within Peano 
Arithmetic. 
\end{fit}

\noindent
[{\bf Query:} Is Fermat's Last Theorem one of these statements?
Wiles has shown it follows from the usual \textsf{ZFC} axioms; but does 
it already follow from the Peano axioms?]

\begin{sit}
Every recursively axiomatizable theory rich enough to include the Peano 
Axioms is incapable of demonstrating its own consistency. 
\end{sit}

Another fundamental discovery of G\"odel was:

\begin{ct}
Every consistent set of  axioms has a model. 
\end{ct}

Together with the first incompleteness theorem, this has been a source of
a wealth of mathematics as well as such paradoxical facts as the
following: it is impossible to unambiguously formalize the distinction
between ``finite'' and ``infinite''.  The ``featherless biped'' definition
of an infinite set as one that can be put into one-to-one correspondence
with a proper subset of itself does not work; neither does the more
natural definition of a finite set as one that can be put into one-to-one
correspondence with $\{0, \dots, n\}$ for some natural number $n$: the
very concept of ``natural number'' cannot be formalized in a way that
makes it clear that our intuitive concept of a natural number is intended.

\textsc{Go\"del} [1940] also gave a partial solution to Hilbert's First
Problem by showing that the Continuum Hypothesis (\textsf{CH}) is
consistent if the usual Zermelo-Fraenkel (\textsf{ZF}) axioms for set
theory are consistent.  He produced a model, known as the Constructible
Universe, of the \textsf{ZF} axioms in which both the Axiom of Choice
(\textsf{AC}) and the \textsf{CH} hold.  Then Cohen showed in
\oldstylenums{1963} that the negations of these axioms are also consistent
with \textsf{ZF}; in particular, \textsf{CH} can fail while \textsf{AC}
holds in a model of \textsf{ZF}. Cohen's technique for producing such
models was generalized by Scott, Solovay, and Shoenfeld and a huge variety
of models of \textsf{ZFC} (\textsf{ZF} plus \textsf{AC}) has been produced
in the years since then, affecting many areas of mathematics. The books by
\textsc{Dales and Woodin, Fremlin, Kunen, Kunen and Vaughan, Monk,} and
\textsc{Rudin} as well as the articles by \textsc{Eklof} and
\textsc{Roitman}, and the articles of Blass reviewed by \textsc{Nyikos}
give some idea of how great a variety of topics these independence results
have been relevant to.

Topology has been affected perhaps more than any other field, and the
following gives a small sample.  Recall the Heine-Borel theorem: {\it
Every open cover of a closed bounded subset of the real line has a finite
subcover}.  The conclusion provides also the definition of {\it A compact
topological space}.  The conclusion of another famous topological theorem,
the Bolzano-Weierstrass theorem, is the basis for a weaker concept:

\begin{definition} 
A topological space is {\it countably compact} if every infinite subset 
has an accumulation point.
\end{definition}

A strengthening of countable compactness, not shared by all compact
spaces, is that of {\it sequential compactness:} every sequence has a
convergent subsequence.

These three concepts agree for all metrizable spaces (those spaces whose
topology is given by a distance function to the non-negative reals that is
symmetric, puts distinct points at a positive distance from each other,
and satisfies the triangle inequality).  Compact metrizable spaces have
lots of other properties not shared by compact topological spaces in
general, so it is perhaps surprising that the question of when a compact
space is metrizable can be very simply settled:

\begin{theorem}[Sneider, 1945]  
A compact space is metrizable if, and only if, it is Hausdorff and has a 
$G_\delta$-diagonal; that is, the diagonal $\{(x, x): x \in X\}$ is a 
countable intersection of open sets. 
\end{theorem}

This theorem was extended to all regular countably compact spaces by 
J.~Chaber in \oldstylenums{1975}.  One might naturally expect these two
theorems to either stand or fall together if ``$G_\delta$-diagonal'' is
weakened to ``small diagonal'':

\begin{definition}
A space has a {\it small diagonal} if, whenever $A$ is an uncountable 
subset of $X \times X$  that is disjoint from the diagonal $\Delta$, there 
is a neighborhood $U$ of $\Delta$ such that $U \setminus A$ is 
uncountable.
\end{definition}

But in fact, this is not the case.  On the one hand, we know that
\textsf{CH} implies every compact Hausdorff space with a small diagonal is
metrizable. We do not know whether \textsf{ZFC} implies this as well; but,
be that as it may, the corresponding statement about regular countably
compact spaces is independent not only of \textsf{ZFC}, but also of
\textsf{CH}. On the one hand, Gary Gruenhage has shown that in a model of
\textsf{CH} constructed by Todd Eisworth and Peter Nyikos, the statement
is true---every countably compact regular space with a small diagonal is
metrizable; on the other hand, Oleg Pavlov has constructed a
counterexample in G\"odel's Constructible Universe, the very model that
originally established the consistency of \textsf{CH}!

\section*{Annotated Bibliography}

\noindent
\textsc{J.~Barwise,} ed. 
{\it Handbook of Mathematical Logic,} 
North-Holland, 1977.

Part B, on set theory, has many consistency and independence results,
including applications to topology.  The article by J.~P.~Burgess in Part
B gives a fine introduction to forcing. The article by Smorynski on
G\"odel's Incompleteness Theorems is one of the very few treatments I have
seen that does not have holes plugged only by hand-waving.  After the
clearest treatment I have ever seen of the Second Incompleteness Theorem,
he even points out, ``In Section 2.1, we have been guilty of cheating in
two places'' and then goes on to make the necessary repairs.  There is
also a significant article by Harrington about an incompleteness in Peano
Arithmetic at the end.

\vskip\baselineskip

\noindent
\textsc{F.~Browder,} ed. 
{\it Mathematical Developments Arising from Hilbert Problems,} Proceedings 
of Symposia in Pure Mathematics XXVIII, American Mathematical Society, 1974.  

Includes a reprint of the English translation of Hilbert's article.  The
article on Hilbert's first problem, by D.~A.~Martin, expounds on the
significance of consistency and independence proofs, and of large cardinal
axioms. There are articles on the second problem by Kreisel and on the
tenth by the co-solvers, Martin Davis (see reference to an article by him
below), Yuri Matijasevic, and Julia Robinson. A quote from their article:
``The consistency of a recursively axiomatizable theory is equivalent to
the assertion that some definite Diophantine equation has no solutions.''

\vskip\baselineskip
   
\noindent
\textsc{H.G.~Dales and W.H.~Woodin,} 
{\it An Introduction to Independence for Analysts,} Cambridge University 
Press, 1987.

An eloquent preface introduces a self-contained treatment of the
set-theoretic independence of a basic problem in functional analysis: 
If $X$ is compact, Hausdorff and infinite, is every homomorphism from 
$C(X, \mathbb{C})$ into any Banach algebra continuous?  Answer: No if
\textsf{CH} for {\it every} such infinite $X$, but it is also consistent
that the answer is Yes for {\it every} such $X$!
  
\vskip\baselineskip
   
\noindent
\textsc{M.~Davis,} 
{\it Hilbert's Tenth Problem is Unsolvable,} 
Amer.~Math.~Monthly {\bf 80} (1973) 233--269.

A spellbinding exposition with complete proofs, not merely of the tenth
problem but about how its solution impacts the foundations of mathematics
in completely unexpected ways.  Included is a very concrete treatment of
G\"odel's First Incompleteness Theorem in terms of Diophantine equations.  
If we ever contact an extraterrestrial intelligence and want to impress it
with what human beings are capable of, this would be the article I'd
recommend to be transmitted to them.

\vskip\baselineskip

\noindent
\textsc{P. Eklof,} 
{\it Whitehead's problem is undecidable,} 
Amer.~Math.~Monthly {\bf 83} (1976) 775--788.

The set-theoretic independence of the problem of whether every  Whitehead 
group is free.

\vskip\baselineskip

\noindent
\textsc{D.H.~Fremlin,} 
{\it Consequences of Martin's Axiom,} Cambridge University Press, 1984.

Includes many applications to topology, measure theory, and algebra of
Martin's Axiom and the negation of \textsf{CH}, as well as of some weaker
axioms which also deny \textsf{CH}.

\vskip\baselineskip

\noindent 
\textsc{K.~G\"odel,}
{\it The Consistency of the Continuum Hypothesis,} 
Ann.~Math.~Studies no.~3, Princeton University Press, 1940.\\ 
{\it What is Cantor's continuum problem?} 
Amer.~Math.~Monthly {\bf 54} (1947) 515--525.

The first article gives the proof of its main results in full; the second
explains, {\it inter alia}, why G\"odel believed the Continuum Hypothesis
to be ``dubious'' in spite of its consistency.

\vskip\baselineskip
  
\noindent 
\textsc{A.~Kanamori and M.~Magidor,} 
{\it The evolution of large cardinal axioms in set theory,} pp.~99--275 
in: {\it Higher Set Theory\/}, G.~H.~Muller and D.~S.~Scott, eds., Lecture 
Notes in Math. no.~669, Springer-Verlag, 1978.

A dramatic article on large cardinal axioms with a wealth of information 
and proofs.

\vskip\baselineskip

\noindent
\textsc{K.~Kunen,} 
{\it Set Theory: An Introduction to Independence Proofs}, North-Holland, 1980.

Together with Burgess's article referenced above, this provides a fine
understanding of how forcing is done and why its results are consistent
with \textsf{ZFC}.

\vskip\baselineskip

\noindent 
\textsc{K.~Kunen and J.~Vaughan,} eds. 
{\it Handbook of Set-Theoretic Topology}, North-Holland, 1984.

Still the most comprehensive single source of information about the
subject.

\vskip\baselineskip   

\noindent 
\textsc{P.~Maddy,} 
{\it Believing the axioms. I\/} and {\it Believing the axioms. II,\/}
J.~Symbolic Logic {\bf 53} (1988) 481--511 and 736--764.

A highly readable pair of articles in which a philosopher looks at
\textsf{CH} and at large cardinal axioms, and reasons for believing or
disbelieving them.

\vskip\baselineskip

\noindent 
\textsc{J.D.~Monk,} 
{\it Cardinal Invariants on Boolean algebras,}
Birkh\"auser Verlag, 1996.

Contains many consistency and independence results.

\vskip\baselineskip

\noindent 
\textsc{P.~Nyikos,}
untitled review, J.~Symbolic Logic {\bf 57} (1992) 763--766.

A review of seven papers authored or co-authored by Andreas Blass, giving
applications of forcing to algebra, analysis, and topology.

\vskip\baselineskip
 
\noindent
\textsc{J.~Roitman,} 
{\it The uses of set theory,} Math.~Intelligencer {\bf 14} (1) (1992) 63--69.

An entertaining and informative article which pointedly omits all
applications to general topology, Boolean algebra, Whitehead groups, and
measure theory, in order to better make the point that set-theoretic
consistency results and specialized set-theoretic techniques are useful in
unexpected places in mathematics.

\vskip\baselineskip

\noindent
\textsc{M.E.~Rudin,} 
{\it Lectures on Set Theoretic Topology,}
American Mathematical Society, 1975.

This booklet made it clear how profoundly general topology was being
remade by set-theoretic consistency and independent results.

\vskip\baselineskip

\noindent
\textsc{S.G. Simpson,} 
{\it Partial Realizations of Hilbert's Program,}
J. Symbolic Logic {\bf 53} (1988) 349--363.

After the eloquent words at the beginning which I quoted, Simpson explains
Hilbert's Program for salvaging the foundations of mathematics, and goes
on to show how, despite G\"odel's negative solution of Hilbert's Second
Problem, a lot can be done in this direction.  In particular, he recounts
how a lot of familiar results in analysis can be proven in a system called
\textsf{WKL}${}_\mathsf{0}$, and anything that can be proven in this 
system is finitistically reducible in the way Hilbert envisioned.  
There is also a nice introduction to the field of reverse mathematics, 
which deals with the general question: given a theorem in mathematics, 
which set existence axioms are required to prove it?

\end{document}